\documentclass[11pt]{article}
\usepackage{enumerate}
\usepackage{amssymb,a4wide,latexsym,makeidx,epsfig,fleqn}
\usepackage{amsthm}
\usepackage{amsmath}
\usepackage{enumerate}
\usepackage{graphicx}
\usepackage{subfigure}
\usepackage{float}
\usepackage{caption}
\usepackage[colorlinks, linkcolor=black, anchorcolor=black, citecolor=blue]{hyperref}%
\newtheorem{theorem}{Theorem}[section]

\newtheorem{conjecture}{Conjecture}[section]

\newtheorem{lemma}[theorem]{Lemma}

\begin{document}
\textwidth 150mm \textheight 225mm
\title{Extension on spectral extrema of gem-free graph with given size\thanks{Supported by the National Natural Science Foundation of China (No. 12271439).}}
\author{{Yuxiang Liu$^{a,b}$, Ligong Wang$^{a,b,}$\footnote{Corresponding author.}}\\
{\small $^{a}$  School of Mathematics and Statistics, Northwestern Polytechnical University,}\\{\small  Xi'an, Shanxi 710129, P.R. China.}\\ {\small $^{b}$ Xi'an-Budapest Joint Research Center for Combinatorics, Northwestern Polytechnical University,}\\{\small  Xi'an, Shanxi 710129, P.R. China.}\\
{\small E-mail: yxliumath@163.com, lgwangmath@163.com}}
\date{}
\maketitle
\begin{center}
\begin{minipage}{135mm}
\vskip 0.3cm
\begin{center}
{\small {\bf Abstract}}
\end{center}
{\small  A graph $G$ is $F$-free if $G$ does not contain $F$ as a subgraph. Let $\mathcal{G}(m, F)$ denote the family of $F$-free graphs with $m$ edges and without isolated vertices. Let $S_{n,k}$ denote the graph obtained by joining every vertex of $K_{k}$ to $n-k$ isolated vertices and $S_{n,k}^{t}$ denote the graph obtained from $S_{n-t,k}$ by attaching $t$ pendant vertices to the maximal degree vertex of $S_{n-t,k}$, respectively. Denote by $H_{n}$ the fan graph obtain from $n-1$-vertex path plus a vertex adjacent to each vertex of the path. Particularly, the graph $H_{5}$ is also known as the gem. Zhang and Wang [Discrete Math. 347(2024)114171] and Yu, Li and Peng [arXiv: 2404. 03423] showed that every gem-free graph $G$ with $m$ edges satisfies $\rho(G)\leq \rho(S_{\frac{m+3}{2},2})$. In this paper, we show that if $G\in \mathcal{G}(m, H_{5})\setminus S_{\frac{m+3}{2},2}$ be a graph of odd size $m\geq23$, then $\rho(G)\leq \rho(S_{\frac{m+5}{2},2}^{2})$, and equality holds if and only if $G\cong S_{\frac{m+5}{2},2}^{2}$.

}
\end{minipage}
\end{center}
\section{Introduction}
Throughout this paper, we consider all graphs are always undirected and simple. We follow the traditional notation and terminology \cite{BoMu1}. Let $G$ be a graph of order $n$ with vertex set $V(G)=\{v_{1}, v_{2},\ldots, v_{n}\}$ and size $m$ with edge set $E(G)$. For a vertex $u\in V(G)$, let $N_{G}(u)$ be the neighborhood set of a vertex $u$, $N_{G}[u]=N_{G}(u)\cup \{u\}$ and $N_{G}^{2}(u)$ be the set of vertices of distance two to $u$ in $G$. In particular, $N_{S}(v)=N(v)\cap S$ and $d_{S}(v)=|N_{S}(v)|$ for a subset $S\subseteq V(G)$. Let $d_{G}(u)=|N_{G}(u)|$ be the degree of a vertex $u$. For the sake of simplicity, we omit all the subscripts if $G$ is clear from the context, for example, $N(u)$, $N[u]$, $N^{2}(u)$ and $d(u)$. For a graph $G$ and a subset $S\subseteq V(G)$, let $G[S]$ be the subgraph of $G$ induced by $S$. For two vertex subsets $S$ and $T$ of $V(G)$ (where $S\cap T$ may not be empty), let $e(S,T)$ denote the number of edges with one endpoint in $S$ and the other in $T$. $e(S,S)$ is simplified by $e(S)$. Given two vertex-disjoint graphs $G$ and $H$, we denote by $G\cup H$ the disjoint union of the two graphs, and by $G\vee H$ the joint graph obtained from $G\cup H$ by joining each vertex of $G$ with each vertex of $H$. The adjacency matrix of a graph $G$ is an $n\times n$ matrix $A(G)$ whose $(i,j)$-entry is $1$ if $v_{i}$ is adjacent to $v_{j}$ and $0$ otherwise. The spectral radius $\rho(G)$ of a graph $G$ is the largest eigenvalue of its adjacency matrix $A(G)$.

Let $P_{n}, C_{n}, K_{1,n-1}$ and $K_{t,n-t}$ be the path of order $n$, the cycle of order $n$, the star graph of order $n$ and the complete bipartite graph with two parts of sizes $t, n-t$, respectively. Let $S_{n}^{k}$ be the graph obtained from $K_{1,n-1}$ by adding $k$ disjoint edges within its independent sets. Let $S_{n,k}$ be the graph obtained by joining every vertex of $K_{k}$ to $n-k$ isolated vertices. Let $S_{n,k}^{t}$ be the graph obtained from $S_{n-t,k}$ by attaching $t$ pendant vertices to the maximal degree vertex of $S_{n-t,k}$. Denote by $H_{n}$ the fan graph obtain from $n-1$-vertex path plus a vertex adjacent to each vertex of the path. Particularly, the graph $H_{5}$ is also known as the gem.

Given a graph $F$, a graph $G$ is $F$-free if it does not contain $F$ as a subgraph. Let $\mathcal{G}(m, F)$ denote the family of $F$-free graphs with $m$ edges and without isolated vertices. In spectral graph theory, Nikiforov \cite{Ni2} posed a spectral Tur\'{a}n type problem which asks to determine the maximum spectral radius of an $F$-free graphs of $n$ vertices, which is known as the Brualdi-Solheid-Tur\'{a}n type problem. In the past decades, this problem has received much attention. For more details, we suggest the reader to see surveys \cite {ChZ,FSi,LiLF,Ni3}, and references therein. In addition, Brualdi and Hoffman \cite{BrH} raised another spectral Tur\'{a}n type problem: What is the maximal spectral radius of an $F$-free graph with given size $m$? We mainly pay attention to this the stability results of Brualdi-Hoffman Tur\'{a}n type problem in this paper. Nosal \cite{Nos} showed that every triangle-free graph $G$ with $m$ edges satisfies $\rho(G)\leq\sqrt{m}$. Lin, Ning and Wu \cite{LinNW} slightly improved the bound to $\rho(G)\leq\sqrt{m-1}$ when $G$ is non-bipartite and triangle-free graph of size $m$. Zhai and Shu \cite{ZhS} further showed the sharp bound to $\rho(G)\leq \rho(SK_{2,\frac{m-1}{2}})$ when $G$ is non-bipartite and triangle-free graph of size $m$ (where $SK_{2,\frac{m-1}{2}}$ be the graph obtained by subdividing an edge of $K_{2,\frac{m-1}{2}}$). Nikiforov showed that every quadrilateral-free graph $G$ with $m$ edges satisfies $\rho(G)\leq\rho(K_{1,m})$. Zhai and Shu \cite{ZhS} further showed the sharp bound to $\rho(G)\leq \rho(S_{m}^{1})$ when $G$ is non-bipartite and quadrilateral-free graph of size $m$. Zhai, Lin and Shu \cite{ZhLS} showed that every $C_{5}/C_{6}$-free graph $G$ with size $m$ satisfies $\rho(G)\leq \rho(S_{\frac{m+3}{2},2})$. Sun, Li and Wei \cite{SLW} improved the bound to $\rho(G)\leq \rho(S_{\frac{m+5}{2},2}^{2})$ when $G\ncong S_{\frac{m+3}{2},2}$ and $C_{5}/C_{6}$-free graph of size $m$. For the Brualdi-Solheid-Tur\'{a}n type problem of fan-free graph of given size $m$, Yu, Li and Peng \cite{YuLP} proposed the following Conjecture.

\noindent\begin{conjecture}\label{co:ch-1.1.}{\rm(}$\cite{YuLP}${\rm)}
For $k\geq2$ and $m$ sufficiently large, if $G$ is an $H_{2k+1}$-free or $H_{2k+2}$-free  with given size $m$, we have $\rho(G)\leq\frac{k-1+\sqrt{4m-k^{2}+1}}{2}$, and the equality holds if and only if $G\cong K_{k}\vee (\frac{m}{k}-\frac{k-1}{2})K_{1}$.
\end{conjecture}

Recently, Li, Zhou and Zou \cite{LiZZ} completely solved Conjecture \ref{co:ch-1.1.}. Zhang and Wang \cite{ZhW} and Yu, Li and Peng \cite{YuLP} solved Conjecture \ref{co:ch-1.1.} for the case $H_{5}$, respectively.

\noindent\begin{theorem}\label{th:ch-1.1.}{\rm(}$\cite{YuLP,ZhW}${\rm)}
Let $G\in \mathcal{G}(m, H_{5})$ be a graph of size $m\geq11$, then $\rho(G)\leq \rho(S_{\frac{m+3}{2},2})$, and equality holds if and only if $G\cong S_{\frac{m+3}{2},2}$.
\end{theorem}

We shall notice that the size of the unique maximal graph in Theorem 1 is odd. Recently, Chen and Yuan showed that the Brualdi-Solheid-Tur\'{a}n type problem of gem-free graph of given even size $m$. Motivated by their result, in what follows, we will further show the sharp bound of spectral radius when $G\in \mathfrak{G}(m, H_{5})$ is $H_{5}$-free graph of odd size $m$ (where $\mathfrak{G}(m, H_{5})=\mathcal{G}(m, H_{5})\setminus \{S_{\frac{m+3}{2},2}\}$).

\noindent\begin{theorem}\label{th:ch-1.2.}
Let $G\in \mathfrak{G}(m, H_{5})$ be a graph of odd size $m\geq23$, then $\rho(G)\leq \rho(S_{\frac{m+5}{2},2}^{2})$, and equality holds if and only if $G\cong S_{\frac{m+5}{2},2}^{2}$.
\end{theorem}

\section{Preliminary}
In this section, we present some preliminary results which used to prove our results.

\noindent\begin{lemma}\label{le:ch-2.1.} {\rm(}$\cite{ZhWF1}${\rm)} Let $u, v$ be two distinct vertices of the connected graph $G$, $\{v_{i}| i=1,2, \ldots, s\}\subseteq N(v)\setminus N(u)$, and $X=(x_{1}, x_{2}, \ldots, x_{n})^{T}$ be the Perron vector of $G$. Let $G^{\prime}=G-\sum_{i=1}^{s}v_{i}v+\sum_{i=1}^{s}v_{i}u$. If $x_{u}\geq x_{v}$, then $\rho(G)<\rho(G^{\prime})$.
\end{lemma}

\noindent\begin{lemma}\label{le:ch-2.2.}{\rm(}$\cite{SLW}${\rm)} Let $m(\geq22)$ and $t$ be two positive integers such that $m\geq t+1$.

(i)$\rho(S_{\frac{m+5}{2},2}^{2})>\frac{1+\sqrt{4m-7}}{2}$ ;

(ii)$\rho(S_{\frac{m+5}{2},2}^{2})>\rho(S_{\frac{m+t+3}{2},2}^{t})$ for even $t(\geq4)$.

\end{lemma}

A cut vertex of a graph is a vertex whose deletion increases the number of components. A graph is called 2-connected, if
it is a connected graph without cut vertices. A block is a maximal 2-connected subgraph of a graph. An end-block is a block
containing at most one cut vertex.

\noindent\begin{lemma}\label{le:ch-2.3.}
Let $G$ be a graph having the maximum spectral radius among $\mathfrak{G}(m, H_{5})$, where $m\geq 7$. Then $G$ is
connected. Furthermore, if $u$ is the extremal vertex of $G$, then there exists no cut vertex in $V(G)\ \{u\}$, and so $d(v)\geq2$ for all $v\in V(G)\ N[u]$.
\end{lemma}

\noindent{\bf Proof.} At first, suppose to the contrary that $G$ is disconnected. Let $G_{1}$ and $G_{2}$ be two connected components of $G$ and $\rho(G)=\rho(G_{1})$. Choose $v_{1}v_{2}\in E(G_{2})$ and assume $v_{3}$ is a vertex in $G_{1}$. Let $G_{3}$ be the graph obtained from $G-v_{1}v_{2}+v_{1}v_{3}$ by deleting all of its isolated vertices. Note that $H_{5}$ is 2-connected and $v_{1}v_{3}$ is a cut edge of $G_{3}$. Hence, $G_{3}$ is $H_{5}$-free and $G_{3}\ncong S_{\frac{m+3}{2},2}$. Clearly, $G_{1}$ is a proper subgraph of $G_{3}$, Therefore, $\rho(G_{3})>\rho(G_{1})=\rho(G)$, which contradicts the choice of $G$. Hence, $G$ is connected.

Suppose that there is at least one cut vertex of $G$ in $V(G)\setminus\{u\}$. Let $B$ be an end-block of $G$ with $u\notin V(B)$ and let $v\in V(B)$ be a cut vertex of $G$. Let $G_{4}=G-\{wv:w\in V(B)\cap N(v)\}+\{wu:w\in V(B)\cap N(v)\}$. Obviously, $u$ is a cut vertex of
$G_{4}$. Then $G_{4}\ncong S_{\frac{m+3}{2},2}$. It is routine to check that $G_{4}$ is $H_{5}$-free. That is, $G_{4}\in\mathfrak{G}(m, H_{5})$. In view of Lemma \ref{le:ch-2.1.}, we have $\rho(G_{4})>\rho(G)$, a contradiction.$\qedsymbol$

\section{Proof of Theorem 1.2.}
In this section, we give the proof of Theorem 1.2, and assume that $\hat{G}$ is a graph in $\mathfrak{G}(m,H_{5})$ such that $\hat{\rho}:=\rho(\hat{G})$ is as large as possible. In view of Lemma \ref{le:ch-2.3.}, we know that $\hat{G}$ is connected. Assume that $\hat{x}$ is the Perron vector of $\hat{G}$. A vertex $\hat{u}$ in $\hat{G}$ is said to be an extremal vertex if $x_{\hat{u}}=max\{x_{v}\mid v\in V(\hat{G})\}$. Notice that $S_{\frac{m+5}{2},2}^{2}$ is
$H_{5}$-free. By Lemma \ref{le:ch-2.2.}, we have $\hat{\rho}>\rho(S_{\frac{m+5}{2},2}^{2})>\frac{1+\sqrt{4m-7}}{2}$, that is ${\hat{\rho}}^{2}-\hat{\rho}>m-2$. Since $\hat{G}$ is $H_{5}$-free, we know that $\hat{G}[N(\hat(u))]$ contains no path of length 3. Then the following result holds immediately.

\noindent\begin{lemma}\label{le:ch-3.1.} Each connected component of $\hat{G}[N(\hat(u))]$ is either a triangle or a star $K_{1,r}$ for  some integer $r\geq0$, where $K_{1,0}$ is a singleton component.
\end{lemma}

For convenience, let $W=V(\hat{G})\backslash N[\hat{u}]$, $W_{0}=\cup_{u\in N_{0}(\hat{u})}N_{W}(u)$ and
$N_{+}(\hat{u})=N(\hat{u})\setminus N_{0}(\hat{u})$, where $N_{0}(\hat{u})$ denotes the set of isolated vertices of $\hat{G}[N(\hat{u})]$.

Since $A(\hat{G})X=\hat{\rho} X$, we have

$$\hat{\rho} x_{\hat{u}}=\sum_{u\in N_{+}(\hat{u})}x_u+\sum_{u\in N_{0}(\hat{u})}x_{u}.$$

Furthermore, we have
\begin{equation}\label{eq:ch-1}
\begin{aligned}
\hat{\rho}^{2}x_{\hat{u}}&=\sum_{u\in N_{+}(\hat{u})}\hat{\rho}x_{u} +\sum_{u\in N_{0}(\hat{u})}\hat{\rho} x_{u}
\\&=|N_{+}(\hat{u})|x_{\hat{u}}+\sum_{u\in N_{+}(\hat{u})}d_{N(\hat{u})}(u)x_{u}+\sum_{w\in N_{W}(u)}d_{W}(u)x_w+|N_{0}(\hat{u})|x_{\hat{u}}
\\&=d(\hat{u})x_{\hat{u}}+\sum_{u\in N_{+}(\hat{u})}d_{N(\hat{u})}(u)x_{u}+\sum_{w\in N^{2}(\hat{u})}d_{N(\hat{u})}(w)x_w.
\end{aligned}
\end{equation}
Therefore,
\begin{equation}\label{eq:ch-2}
\begin{aligned}
(\hat{\rho}^{2}-\hat{\rho})x_{\hat{u}}&=d(\hat{u})x_{\hat{u}}+\sum_{v\in N_{+}(\hat{u})}(d_{N(\hat{u})}(v)-1)x_{v}+\sum_{w\in N^{2}(\hat{u})}d_{N(\hat{u})}(w)x_w-\sum_{v\in N_{0}(\hat{u})}x_{v}
\\&\leq |N(\hat{u})|x_{\hat{u}}+\sum_{u\in N_{+}(\hat{u})}(d_{N(\hat{u})}(u)-1)x_{u}+e(N(\hat{u}),W)x_{\hat{u}}-\sum_{u\in N_{0}(\hat{u})}x_{u}.
\end{aligned}
\end{equation}

Note that $S_{\frac{m+5}{2},2}^{2}$ is $H_{5}$-free, we have $\rho\geq\rho(S_{\frac{m+5}{2},2}^{2})> \frac{1+\sqrt{4m-7}}{2}>5.1$ for $m\geq23$. Thus,

\begin{equation}\label{eq:ch-3}
(\hat{\rho}^{2}-\hat{\rho})x_{\hat{u}}>(m-2)x_{\hat{u}}=(|N(\hat{u})|+e(N_{+}(\hat{u}))+e(N(\hat{u}),W)+e(W)-2)x_{\hat{u}}. \end{equation}

Combining with \eqref{eq:ch-2} and \eqref{eq:ch-3}, we get
\begin{equation}\label{eq:ch-4}
\begin{aligned}
e(W)&<\sum_{u\in N_{+}(\hat{u})}(d_{N(\hat{u})}(u)-1)\frac{x_{u}}{x_{\hat{u}}}-e(N_{+}(\hat{u}))+2-\sum_{u\in N_{0}(\hat{u})}\frac{x_{u}}{x_{\hat{u}}}\\
&\leq e(N_{+}(\hat{u}))-|N_{+}(\hat{u})|+2-\sum_{u\in N_{0}(\hat{u})}\frac{x_{u}}{x_{\hat{u}}}
\end{aligned}
\end{equation}

\noindent\begin{lemma}\label{le:ch-3.2.} $e(N_{+}(\hat{u}))\geq4$.
\end{lemma}
\noindent{\bf Proof.} Suppose to the contrary that $e(N_{+}(\hat{u}))\leq3$. By \eqref{eq:ch-1}, we have
$$
\begin{aligned}
\hat{\rho}^{2}x_{\hat{u}}&=d(\hat{u})x_{\hat{u}}+\sum_{u\in N_{+}(\hat{u})}d_{N(\hat{u})}(u)x_{u}+\sum_{w\in N^{2}(\hat{u})}d_{N(\hat{u})}(w)x_w\\
&\leq\left(|N(\hat{u})|+2e(N_{+}(\hat{u}))+e(N(\hat{u}),W) \right)x_{\hat{u}}\\
&=(m-e(W)+e(N_{+}(\hat{u})))x_{\hat{u}}\\
&\leq (m+3)x_{\hat{u}},
\end{aligned}
$$
which yields $\hat{\rho}\leq \sqrt{m+3}<\frac{1+\sqrt{4m-7}}{2}$ for $m\geq23$, a contradiction.$\qedsymbol$

\noindent\begin{lemma}\label{le:ch-3.3.} If there exists a subset $S\subseteq N_{+}(\hat{u})$ such that $x_{u}<(1-\beta)x_{\hat{u}}$ and $d_{N(\hat{u})}(u)\geq2$ for any vertex $u\in S$, where $0<\beta<1$, then
$$0< e(N_{+}(\hat{u}))-|N_{+}(\hat{u})|+2-\beta\sum_{u\in S}(d_{N(\hat{u})}(u)-1).$$
\end{lemma}
\noindent{\bf Proof.}
$$
\begin{aligned}
\sum\limits_{u\in S}(d_{N(\hat{u})}(u)-1)x_{u}&=\sum\limits_{v\in N_{+}(\hat{u})\setminus S}(d_{N(\hat{u})}(v)-1)x_{v}+\sum\limits_{u\in S}(d_{N(\hat{u})}(u)-1)x_{u}\\
&<\sum\limits_{v\in N_{+}(\hat{u})\setminus S}(d_{N(\hat{u})}(v)-1)x_{v}+(1-\beta)\sum\limits_{u\in S}(d_{N(\hat{u})}(u)-1)x_{\hat{u}}\\
&\leq \sum\limits_{v\in N_{+}(\hat{u})}(d_{N(\hat{u})}(v)-1)x_{\hat{u}}-\beta\sum\limits_{u\in S}(d_{N(\hat{u})}(u)-1)x_{\hat{u}}.
\end{aligned}
$$
Together \eqref{eq:ch-1}, we obtain

$$
\begin{aligned}
(m-2)x_{\hat{u}}&<|N(\hat{u})|x_{\hat{u}}+\sum\limits_{u\in N_{+}(\hat{u})}(d_{N(\hat{u})}(u)-1)x_{u}+e(N(\hat{u}),W)x_{\hat{u}}\\
&<|N(\hat{u})|x_{\hat{u}}+\sum\limits_{v\in N_{+}(\hat{u})}(d_{N(\hat{u})}(v)-1)x_{\hat{u}}-\beta\sum\limits_{u\in S}(d_{N(\hat{u})}(u)-1)x_{\hat{u}}+e(N(\hat{u}),W)x_{\hat{u}}\\
&\leq \left(|N(\hat{u})|+2e(N(\hat{u}))-|N_{+}(\hat{u})|-\beta\sum\limits_{u\in S}(d_{N(\hat{u})}(u)-1)+e(N(\hat{u}),W)\right)x_{\hat{u}}\\
&=\left( m-e(W)+e(N(\hat{u}))-|N_{+}(\hat{u})|-\beta\sum\limits_{u\in S}(d_{N(\hat{u})}(u)-1)\right)x_{\hat{u}},
\end{aligned}
$$
which yields that $0\leq e(W)< e(N_{+}(\hat{u}))-|N_{+}(\hat{u})|+2-\beta\sum\limits_{u\in S}(d_{N(\hat{u})}(u)-1)$.$\qedsymbol$

\noindent\begin{lemma}\label{le:ch-3.4.} Suppose that there are three distinct vertices $u_{1},u_{2},u_{3}\in N_{+}(\hat{u})$ forming a $K_{3}$. Then the following statement hold.

(i) If $\mathop{\cup}\limits_{i=1}\limits^{3}N_{W}(u_{i})=\emptyset$, then $x_{u_{i}}<\frac{10}{31}x_{\hat{u}}$;

(ii) If $\mathop{\cup}\limits_{i=1}\limits^{3}N_{W}(u_{i})\neq \emptyset$ and $x_{u_{1}}\leq x_{u_{2}}\leq x_{u_{3}}$, then $x_{u_{i}}<\frac{20}{41}x_{\hat{u}}$ for each $i\in\{1,2\}$ and $N(w)\cap \{u_{1},u_{2},u_{3}\}=\{u_{3}\}$ for any vertex $w\in\mathop{\cup}\limits_{i=1}\limits^{3}N_{W}(u_{i})$.
\end{lemma}

\noindent{\bf Proof.} (i) If $\mathop{\cup}\limits_{i=1}\limits^{3}N_{W}(u_{i})=\emptyset$, then $$x_{u_{1}}= x_{u_{2}}= x_{u_{3}}, \quad \hat{\rho}x_{u_{1}}=x_{u_{2}}+ x_{u_{3}}+x_{\hat{u}}.$$

Together with $\hat{\rho}>\frac{1+\sqrt{4m-7}}{2}>5.1$ for $m\geq23$, we obtain $x_{u_{1}}< \frac{1}{\hat{\rho}-2}x_{\hat{u}}<\frac{10}{31}x_{\hat{u}}$.

(ii) Suppose that $\mathop{\cup}\limits_{i=1}\limits^{3}N_{W}(u_{i})\neq \emptyset$. For any $w\in\mathop{\cup}\limits_{i=1}\limits^{3}N_{W}(u_{i})$, since $\hat{G}$ is $H_{5}$-free, we have $|N(w)\cap \{u_{1},u_{2},u_{3}\}|=1$. Otherwise suppose $w\sim u_{1}$ and $w\sim u_{2}$. Hence, $\hat{G}[u_{1},\hat{u},u_{3},u_{2},w]$ is the $H_{5}$, a contradiction. We will prove that $N(w)\cap \{u_{1},u_{2},u_{3}\}=\{u_{3}\}$. Without loss of generality, suppose to the contrary that $N(w)\cap \{u_{1},u_{2},u_{3}\}=\{u_{1}\}$. Let $G^{\prime}=\hat{G}-wu_{1}+wu_{3}$. By \eqref{eq:ch-4}, we have $e(W)\leq1$. Then $G^{\prime}$ is still $H_{5}$-free with size $m$. Since $x_{u_{3}}\geq x_{u_{1}}$, combining Lemma \ref{le:ch-2.1.}, we have $\rho(G^{\prime})>\hat{\rho}$, a contradiction. Thus, $N(w)\cap \{u_{1},u_{2},u_{3}\}=\{u_{3}\}$. Furthermore, $$\hat{\rho}x_{u_{1}}=x_{u_{2}}+x_{u_{3}}+x_{\hat{u}},$$

Together with $\hat{\rho}>\frac{1+\sqrt{4m-7}}{2}>5.1$ for $m\geq23$, we obtain $x_{u_{1}}< \frac{2}{\hat{\rho}-1}x_{\hat{u}}<\frac{20}{41}x_{\hat{u}}$.$\qedsymbol$

\noindent\begin{lemma}\label{le:ch-3.5.} If there are three distinct vertices $u_{1},u_{2},u_{3}\in N_{+}(\hat{u})$ forming a $K_{3}$, then $$\mathop{\cup}\limits_{i=1}\limits^{3}N_{W}(u_{i})\neq \emptyset.$$
\end{lemma}

\noindent{\bf Proof.} Suppose to the contrary that $\mathop{\cup}\limits_{i=1}\limits^{3}N_{W}(u_{i})=\emptyset$. By Lemmas \ref{le:ch-3.3.} and \ref{le:ch-3.4.}(i), we obtain $$0< e(N_{+}(\hat{u}))-|N_{+}(\hat{u})|+2-\frac{21}{31}\sum_{u\in S}(d_{N(\hat{u})}(u)-1)=2-\frac{63}{31}<0,$$ a contradiction.$\qedsymbol$

\noindent\begin{lemma}\label{le:ch-3.6.}

(i) $\hat{G}[N_{+}(\hat{u})]$ contains at most one star-component;

(ii) If $\hat{G}[N_{+}(\hat{u})]$ isomorphic to a star, then $e(W)=0$.
\end{lemma}

\noindent{\bf Proof.} (i) Suppose to the contrary that $\hat{G}[N_{+}(\hat{u})]$ contains at least two star-component. By \eqref{eq:ch-4}, we have $$0\leq e(W)<e(N_{+}(\hat{u}))-|N_{+}(\hat{u})|+2-\sum_{u\in N_{0}(\hat{u})}\frac{x_{u}}{x_{\hat{u}}}=-\sum_{u\in N_{0}(\hat{u})}\frac{x_{u}}{x_{\hat{u}}}\leq0,$$ a contradiction. Thus, (i) holds.

(ii) If $\hat{G}[N_{+}(\hat{u})]$ isomorphic to a star, then we have  $e(W)<e(N_{+}(\hat{u}))-|N_{+}(\hat{u})|+2-\sum_{u\in N_{0}(\hat{u})}\frac{x_{u}}{x_{\hat{u}}}=1-\sum_{u\in N_{0}(\hat{u})}\frac{x_{u}}{x_{\hat{u}}}\leq1$ from \eqref{eq:ch-3},, which yields that $e(W)=0$. Thus, (ii) holds.$\qedsymbol$

\noindent\begin{lemma}\label{le:ch-3.7.} $W=\emptyset$.
\end{lemma}

\noindent{\bf Proof.} Suppose to the contrary that $W\neq \emptyset$. We give the following two claims.

\noindent{\bf Claim 3.1.} $\hat{G}[N_{+}(\hat{u})]$ is isomorphic to a star $K_{1,r}$ where $r\geq4$.

\noindent{\bf Proof.} We first show that
\begin{equation}\label{eq:ch-5}
e(N_{+}(\hat{u}))-|N_{+}(\hat{u})|\leq -1.
\end{equation}
Suppose to the contrary that $e(N_{+}(\hat{u}))-|N_{+}(\hat{u})|=0$. Then we obtain each component of $\hat{G}[N_{+}(\hat{u})]$ is isomorphic to a $K_{3}$ from Lemma \ref{le:ch-3.1.}. Furthermore, we assert that $\hat{G}[N_{+}(\hat{u})]$ is $2K_{3}$-free. Otherwise, let $\hat{G}[\{u_{1},u_{2},u_{3}\}]$ and $\hat{G}[\{v_{1},v_{2},v_{3}\}]$  be two triangles. By Lemmas \ref{le:ch-3.3.} and \ref{le:ch-3.4.} (ii), we have
$$0<e(N_{+}(\hat{u}))-|N_{+}(\hat{u})|+2-\frac{21}{41}
\sum_{i=1}^{2}(d_{N(\hat{u})}(u_{i})-1)-\frac{21}{41}\sum_{i=1}^{2}(d_{N(\hat{u})}(_v{i})-1)
=2-\frac{84}{41}<0,$$
a contradiction. Thus, $\hat{G}[N_{+}(\hat{u})]\cong K_{3}$. Let $N_{+}(\hat{u})=\{z_{1},z_{2},z_{3}\}$ with $z_{1}\leq z_{2}\leq z_{3}$. By Lemmas \ref{le:ch-3.4.} and \ref{le:ch-3.5.}, we have $N_{\hat{G}}(w)\cap\{z_{1},z_{2},z_{3}\}=\{z_{3}\}$. Let $G^{\star}=\hat{G}-wz_{3}+w\hat{u}$. Clearly, $G^{\star}$ is $H_{5}$-free. By Lemma \ref{le:ch-2.1.}, we have $\rho(G^{\star})>\hat{\rho}$, a contradiction. Thus, \eqref{eq:ch-5} holds.

We next show that $\hat{G}[N_{+}(\hat{u})]$ is $K_{3}$-free. Otherwise, let $\hat{G}[\{z_{1},z_{2},z_{3}\}]$ is a triangle. By Lemmas \ref{le:ch-3.3.} and \ref{le:ch-3.4.}  (ii) ,
we have $$0<e(N_{+}(\hat{u}))-|N_{+}(\hat{u})|+2
-\frac{21}{41}\sum\limits_{i=1}\limits^{2}(d_{N(\hat{u})}(u_{i})-1)=1-\frac{42}{41}
<0,$$ a contradiction. Thus, $\hat{G}[N_{+}(\hat{u})]$ is $K_{3}$-free. By Lemmas \ref{le:ch-3.1.}, \ref{le:ch-3.2.} and  \ref{le:ch-3.6.} (i), we obtain that $\hat{G}[N_{+}(\hat{u})]\cong K_{1,r}$ for $r\geq4$ and hence $\hat{G}[N(\hat{u})]\cong K_{1,r}\cup tK_{1}$ for $r\geq4$ and $t\geq0$. It follows from Lemma \ref{le:ch-3.6.} (ii) that $e(W)=0$. Furthermore, by \eqref{eq:ch-4}, we have $$0=e(W)<e(N_{+}(\hat{u}))-|N_{+}(\hat{u})|+2-\sum\limits_{u\in N_{0}(\hat{u})}\frac{x_{u}}{x_{\hat{u}}}=1-\sum\limits_{u\in N_{0}(\hat{u})}\frac{x_{u}}{x_{\hat{u}}}$$
and hence $$\sum\limits_{u\in N_{0}(\hat{u})}\frac{x_{u}}{x_{\hat{u}}}<1.$$

Let $V(K_{1,r})=\{v, v_{1}, v_{2},\cdots, v_{r}\}$, where $v$ is the center vertex of the star $K_{1,r}$ and $V(tK_{1})=\{u_{1},u_{2},\cdots, u_{t}\}$. We give the following Claim 3.2.

\noindent{\bf Claim 3.2.} $w\nsim v$ and $u_{j}\nsim w$ for each vertex $w\in W$ and each $j\in\{1,2,\cdots, t\}$.

\noindent{\bf Proof.} Suppose $w\sim v$. We claim that the vertex $w\nsim v_{i}$ for any $i\in\{1,2,\cdots, r\}$. Otherwise, assume that $w\sim v_{1}$. Then $\hat{G}[\{v,\hat{u}, v_{1}, v_{2},w\}]$ is an $H_{5}$, a contradiction. Hence $w\nsim v_{i}$ for any $i\in\{1,2,\cdots, r\}$. Then let $G^{\circ}=\hat{G}-wv+w\hat{u}$. Combining the fact $e(W)=0$, we have $G^{\circ}$ is $H_{5}$-free. By Lemma \ref{le:ch-2.1.}, we have $\rho(G^{\circ})>\hat{\rho}$, a contradiction. Hence, $w\nsim v$.

Note that $\hat{\rho}x_{v}=\sum\limits_{i=1}^{r}x_{v_{i}}+x_{\hat{u}}$ and $\hat{\rho}x_{w}\leq \sum\limits_{i=1}^{r}x_{v_{i}}+\sum\limits_{j=1}^{t}x_{u_{j}}<\sum\limits_{i=1}^{r}x_{v_{i}}+x_{\hat{u}}$ due to $\sum\limits_{u\in N_{0}(\hat{u})}\frac{x_{u}}{x_{\hat{u}}}<1$. Hence $x_{w}< x_{v}$. If $u_{j}\sim w$, then let $G_{5}=\hat{G}-u_{j}w+u_{j}v$. It is checked that $G_{5}$ is still an $H_{5}$-free graph with size $m$. By Lemma \ref{le:ch-2.1.}, we have $\rho(G_{5})>\hat{\rho}$, a contradiction. Hence $u_{j}\nsim w$ for each $j\in\{1,2,\cdots, t\}$.

Let $W=\{w_{1},w_{2},\cdots, w_{c}\}$. From above, we know that $N_{\hat{G}}(w_{i})\subseteq \{v_{1}, v_{2},\cdots, v_{r}\}$, set $d_{i}=d_{\hat{G}}(w_{i})$.  By Lemma \ref{le:ch-2.3.}, $d_{i}\geq2$ for $1\leq i\leq c$. Next we will modify the structure of $\hat{G}$ to obtain a contradiction.
\\
\rule{15cm}{0.3mm}\\
\textbf{Construction 1}:\\
\textbf{Input} Graph $G_{0}=\hat{G}$ and $i=1$.\\
\textbf{Output} Graph $G_{c}$.\\
\textbf{Step 1}  \textbf{While} $d_{i}$ even, \textbf{do}\\
$G_{i}=G_{i-1}\cup \frac{d_{i}}{2}K_{1}$, and $V_{i}=\{k_{1}^{i},k_{2}^{i},\cdots, k_{\frac{d_{i}}{2}}^{i}\}$.\\
\textbf{Otherwise} $d_{i}$ odd, \textbf{do} \\
$G_{i}=G_{i-1}\cup (\frac{d_{i}-1}{2}+1)K_{1}$, and $V_{i}=\{k_{0}^{i},k_{1}^{i},\cdots, k_{\frac{d_{i}-1}{2}}^{i}\}$.\\
\textbf{Step 2}  \textbf{While} $i<c$, \textbf{do} \\
$i=i+1$ and back to \textbf{Step 1}.\\
\textbf{Otherwise}, stop the construction.\\
\rule{15cm}{0.2mm}\\
\textbf{Construction 2}:\\
\textbf{Input} Graph $G_{0}^{\prime}=G_{c}$ and $i=1$.\\
\textbf{Output} Graph $G_{c}^{\prime}$.\\
\textbf{Step 1}  \textbf{While} $d_{i}$ even, \textbf{do}\\
$G_{i}^{\prime}=G_{i-1}^{\prime}-\sum\limits_{v_{j}\in N_{\hat{G}}}w_{i}v_{j}+\sum\limits_{j=1}\limits^{\frac{d_{i}}{2}}
(\hat{u}k_{j}^{i}+vk_{j}^{i})$.\\
\textbf{Otherwise} $d_{i}$ odd, \textbf{do} \\
$G_{i}^{\prime}=G_{i-1}^{\prime}-\sum\limits_{v_{j}\in N_{\hat{G}}}w_{i}v_{j}+\sum\limits_{j=1}\limits^{\frac{d_{i}-1}{2}}
(\hat{u}k_{j}^{i}+vk_{j}^{i})+\hat{u}k_{0}^{i}$.\\
\textbf{Step 2}  \textbf{While} $i<c$, \textbf{do} \\
$i=i+1$ and back to \textbf{Step 1}.\\
\textbf{Otherwise}, stop the construction.\\
\rule{15cm}{0.3mm}\\

By Construction 1, we have $\rho(\hat{G})=\rho(G_{c})$. Let $G_{c}^{\ast}=G_{c}^{\prime}-\{w_{1},w_{2},\cdots, w_{c}\}$. Then $G_{c}^{\ast}\cong S_{\frac{m+t+3}{2},2}^{t}$ for $t\geq1$. Hence $G_{c}^{\ast}$ is a connected $H_{5}$-free graph with size $m$ and $\rho(G_{c}^{\ast})=\rho(G_{c}^{\prime})$. To obtain a contradiction, we will show $\rho(G_{c}^{\ast})>\rho(\hat{G})$.

Let $Y=(X^{T},0)^{T}$ and $Z$ be a non-negative eigenvector of $G_{c}$ and $G_{c}^{\prime}$ corresponding to $\rho(G_{c})$ and $\rho(G_{c}^{\prime})$. Since $k_{j}^{i}$ is an isolated vertex in $G_{c}$ and $w_{i}$ is an isolated vertex in $G_{c}^{\prime}$ for any $1\leq i\leq c$ and $1\leq j\leq d_{i}$. It follows that $y_{k_{j}^{i}}=z_{w_{i}}=0$ for $0\leq j\leq d_{i}$ and $1\leq i\leq c$. Observe that $N_{G_{c}^{\prime}}(k_{j}^{i})=N_{G_{c}^{\prime}}(v_{l})$ for $j\geq1$ and $1\leq l\leq r$. Hence, $z_{k_{j}^{i}}=z_{v_{l}}$ for $j\geq1$ and $1\leq l\leq r$. Combining the fact $x_{w_{i}}<x_{v}\leq x_{\hat{u}}$, we have $x_{\hat{u}}z_{k_{j}^{i}}-x_{w_{i}}z_{v_{l}}=(x_{\hat{u}}-x_{w_{i}})z_{v_{l}}>0$,
$x_{v}z_{k_{j}^{i}}-x_{w_{i}}z_{v_{l}}=(x_{v}-x_{w_{i}})z_{v_{l}}>0$ for $1\leq i\leq c$ and $1\leq j\leq d_{i}$.

If $d_{i}$ is even, then let
$$
\begin{aligned}
f_{w_{i}}&=\sum\limits_{j=1}\limits^{\frac{d_{i}}{2}}(z_{\hat{u}}y_{k_{j}^{i}}+y_{\hat{u}}z_{k_{j}^{i}}+z_{v}y_{k_{j}^{i}}+y_{v}z_{k_{j}^{i}})-\sum\limits_{j=1}\limits^{d_{i}}(z_{w_{i}}y_{v_{j}}+y_{w_{i}}z_{v_{j}})\\
&=\sum\limits_{j=1}\limits^{\frac{d_{i}}{2}}(y_{\hat{u}}z_{k_{j}^{i}}+y_{v}z_{k_{j}^{i}})-\sum\limits_{j=1}\limits^{d_{i}}y_{w_{i}}z_{v_{j}}\\
&=\sum\limits_{j=1}\limits^{\frac{d_{i}}{2}}(y_{\hat{u}}z_{k_{j}^{i}}-y_{w_{i}}z_{v_{j}})+\sum\limits_{j=1}\limits^{\frac{d_{i}}{2}}(y_{v}z_{k_{j}^{i}}-y_{w_{i}}z_{v_{j+\frac{d_{i}}{2}}})\\
&=\sum\limits_{j=1}\limits^{\frac{d_{i}}{2}}(x_{\hat{u}}z_{k_{j}^{i}}-x_{w_{i}}z_{v_{j}})+\sum\limits_{j=1}\limits^{\frac{d_{i}}{2}}(x_{v}z_{k_{j}^{i}}-x_{w_{i}}z_{v_{j+\frac{d_{i}}{2}}})>0.
\end{aligned}
$$

If $d_{i}$ is odd, then let
$$
\begin{aligned}
f_{w_{i}}&=\sum\limits_{j=1}\limits^{\frac{d_{i}-1}{2}}(z_{\hat{u}}y_{k_{j}^{i}}+y_{\hat{u}}z_{k_{j}^{i}}+z_{v}y_{k_{j}^{i}}+y_{v}z_{k_{j}^{i}})+(z_{\hat{u}}y_{k_{0}^{i}}+y_{\hat{u}}z_{k_{0}^{i}})-\sum\limits_{j=1}\limits^{d_{i}}(z_{w_{i}}y_{v_{j}}+y_{w_{i}}z_{v_{j}})\\
&=\sum\limits_{j=1}\limits^{\frac{d_{i}-1}{2}}(y_{\hat{u}}z_{k_{j}^{i}}+y_{v}z_{k_{j}^{i}})+y_{\hat{u}}z_{k_{0}^{i}}-\sum\limits_{j=1}\limits^{d_{i}}y_{w_{i}}z_{v_{j}}\\
&=\sum\limits_{j=1}\limits^{\frac{d_{i}-1}{2}}(y_{\hat{u}}z_{k_{j}^{i}}-y_{w_{i}}z_{v_{j}})+\sum\limits_{j=1}\limits^{\frac{d_{i}-1}{2}}(y_{v}z_{k_{j}^{i}}-y_{w_{i}}z_{v_{j+\frac{d_{i}}{2}}})+(y_{\hat{u}}z_{k_{0}^{i}}-y_{w_{i}}z_{v_{d_{i}}})\\
&=\sum\limits_{j=1}\limits^{\frac{d_{i}-1}{2}}(x_{\hat{u}}z_{k_{j}^{i}}-x_{w_{i}}z_{v_{j}})+\sum\limits_{j=1}\limits^{\frac{d_{i}-1}{2}}(x_{v}z_{k_{j}^{i}}-x_{w_{i}}z_{v_{j+\frac{d_{i}}{2}}})+(x_{\hat{u}}z_{k_{0}^{i}}-x_{w_{i}}z_{v_{d_{i}}})>0.
\end{aligned}
$$

Thus we have $$(\rho(G_{c}^{\prime})-\rho(G_{c}))Y^{T}Z=Y^{T}(A(G_{c}^{\prime})-A(G_{c}))Z=\sum_{i=1}^{c}f_{w_{i}}>0,$$

and hence $\rho(G_{c}^{\prime})>\rho(G_{c})$. By Lemma \ref{le:ch-2.2.} (ii), we obtain that $\rho(S_{\frac{m+5}{2},2}^{2})>\rho(G_{c}^{\ast})=\rho(G_{c}^{\prime})>\rho(G_{c})=\rho(\hat{G})\geq \rho(S_{\frac{m+5}{2},2}^{2})$, a contradiction.

By Lemmas \ref{le:ch-3.6.} (ii) and \ref{le:ch-3.7.}, we have that $\hat{G}\cong S_{\frac{m+t+3}{2},2}^{t}$ and $t$ is even for odd size $m$. By Lemma \ref{le:ch-2.2.} (ii), we have $\hat{G}\cong S_{\frac{m+5}{2},2}^{2}$. This completes the proof. $\blacksquare$

\section*{Data availability}

No data was used for the research described in the article.

\section*{Declaration of competing interest}

The authors declare that they have no conflict of interest.


\begin{thebibliography}{99}
\bibitem{BoMu1} J.A. Bondy, U.S.R. Murty, Graph Theory, in: Graduate Texts in Mathematics, Vol. 244, Springer, New York, 2008.
\bibitem{BrH} R.A. Brualdi, A.J. Hoffman, On the spectral radius of $(0,1)$ matrices, Linear Algebra Appl. 65 (1985)
133--146.
\bibitem{ChY} F. Chen, X.Y. Yuan, Brualdi--Hoffman--Tur$\acute{a}$n problem of the gem, arXiv: 2411.08345v1.
\bibitem{ChZ} M.Z. Chen, X.-D Zhang, Some new results and problems in spectral extremal graph theory, J. Anhui Univ. Nat. Sci. 42 (2018) 12--25 (in Chinese).
\bibitem{FSi} Z. F\"{u}redi, M. Simonovits, The history of degenerate (bipartite) extremal graph problems, Erd\H{o}s centennial, Bolyai Soc. Math. Stud. 25 (2013) 169--264.
\bibitem{LiZZ}  S. C. Li, S. S. Zhao, L. T. Zou, Spectral extrema of graphs with fixed size: forbidden
fan graph, friendship graph or theta graph, arXiv: 2409. 15918.
\bibitem{LiLF} Y.T. Li, W.J. Liu, L.H. Feng, A survey on spectral conditions for some extremal graph problems, Adv. Math. (China) 51 (2022) 193--258.
\bibitem{LinNW} H.Q. Lin, B. Ning, B.Y.D.R. Wu, Eigenvalues and triangles in graphs, Comb. Probab. Comput. 30 (2) (2021) 258--270.
\bibitem{Ni6} V. Nikiforov, Walks and the spectral radius of graphs, Linear Algebra Appl. 418 (2006) 257--268.
\bibitem{Ni2} V. Nikiforov, The spectral radius of graphs without paths and cycles of specified length, Linear Algebra Appl. 432 (2010) 2243--2256.
\bibitem{Ni3} V. Nikiforov, Some new results in extremal graph theory. Lond. Math. Soc. Lect. Note Ser. 392 (2011) 141--181.
\bibitem{Nos} E. Nosal, Eigenvalues of Graphs, Master's Thesis, University of Calgary, 1970.
\bibitem{SLW} W.T. Sun, S.C. Li, W. Wei, Extensions on spectral extrema of $C_{5}/C_{6}$-free graphs with given size, Discrete Math. 346 (2023) 113591.
\bibitem{YuLP} L. J. Yu, Y. T. Li, Y. J. Peng, Spectral extremal graphs for fan graphs, arXiv: 2404.03423.

\bibitem{ZhLS} M.Q. Zhai, H.Q. Lin, J.L. Shu, Spectral extrema of graphs with fixed size: Cycles and complete bipartite graphs, European J. Combin. 95 (2021) 103322.

\bibitem{ZhS} M.Q. Zhai, J.L. Shu,  A spectral version of Mantel's theorem, Discrete Math. 345 (2022) 112630.

\bibitem{ZhWF1} M.Q. Zhai, B. Wang, Proof of a conjecture on the spectral radius of $C_{4}$-free graphs, Linear Algebra Appl. 437 (2012) 1641--1647.
\bibitem{ZhW} Y.T. Zhang, L.G. Wang, On the spectral radius of graphs without a gem, Discrete Math. 347 (2024) 114171.
\end{thebibliography}
\end{document}